\documentclass[11pt,leqno]{article}

\usepackage{amsmath,amsthm}
\usepackage {latexsym}
\usepackage{amssymb}

\newcommand{\supp}{\mbox{\rm supp}}
\newcommand{\sign}{\mbox{\rm sign}}

\newcommand{\bear}{\begin{eqnarray}}
\newcommand{\eear}{\end{eqnarray}}
\newcommand{\eps}{{\varepsilon}}
\newcommand{\R}{{\mathbb R}}

\newcommand{\Z}{{\mathbb Z}}

\newcommand{\Compl}{{\mathbb C}}

\newcommand{\les}{\lesssim}

\def\Lap{\Delta}
\def\nn{\nonumber}

\newtheorem{theorem}{Theorem}
\newtheorem{lemma}[theorem]{Lemma}
\newtheorem{definition}[theorem]{Definition}

\newtheorem{corollary}[theorem]{Corollary}
\newtheorem{prop}[theorem]{Proposition}
\newtheorem{proposition}[theorem]{Proposition}

\theoremstyle{remark}

\def\Pac{P_{a.c.}}
\def\la{\langle}
\def\ra{\rangle}
\def\norm[#1][#2]{\|#1\|_{#2}}
\def\bignorm[#1][#2]{\Big\|#1\Big\|_{#2}}

\begin{document}

\title{Dispersive bounds for the three-dimensional Schr\"odinger equation
 with almost critical potentials}
\date{September 15, 2004}

\author{Michael\ Goldberg}

\maketitle

\begin{abstract}
We prove a dispersive estimate for the time-independent Schr\"odinger operator
$H = -\Lap + V$ in three dimensions.  The potential $V(x)$ is assumed to lie 
in the intersection $L^p(\R^3)\cap L^q(\R^3)$, $p < \frac32 < q$, and also
to satisfy a generic zero-energy spectral condition.  This class, which
includes potentials that have pointwise decay $|V(x)| \le C(1+|x|)^{-2-\eps}$,
is nearly critical with respect to the natural scaling of the Laplacian.
No additional regularity, decay, or positivity of $V$ is assumed.
\end{abstract}

\section{Introduction}
The propogator $e^{-it\Lap}$ of the free Schr\"odinger equation in $\R^3$
may be represented as a convolution operator with kernel $(4\pi it)^{-3/2}
e^{-i(|x|^2/4t)}$.  From this formula it is clear that the free evolution
satisfies the dispersive bound $\norm[e^{-it\Lap}][1\to\infty] \le 
(4\pi|t|)^{-3/2}$ at all times $t \not= 0$.  In this paper we
consider the perturbed Hamiltonian $H = -\Lap + V$ and seek to prove
similar estimates on the time evolution operator $e^{itH}P_{ac}(H)$.
The projection onto the absolutely continuous spectrum of $H$, denoted here
by $P_{ac}(H)$, is needed to eliminate bound states which do not decay over any
length of time.  Our goal is to avoid placing excessive restrictions on the
regularity, positivity, and decay of the potential $V = V(x)$.  To that end
we formulate the following theorem.

\begin{theorem}
\label{thm:dispersive}
Let $V \in L^p(\R^3) \cap L^q(\R^3)$, $p < \frac32 < q$.
Assume also that zero is neither an eigenvalue nor a resonance of 
$H=-\Lap+V$. Then
\begin{equation} \label{eq:dispersive}
\big\|e^{itH} P_{ac}(H)\big\|_{1\to\infty} \les |t|^{-\frac32}. 
\end{equation}
\end{theorem}

A precise definition of resonances is given in section~\ref{sec:lowenergy}.
With this assumption the spectrum is known to be purely absolutely continuous
on $[0,\infty)$, see \cite{GS2} for details.  We remark that if the 
zero--energy hypothesis is not satisfied, a dispersive estimate still holds
for $e^{itH}P_{[a,\infty)}(H)$ for any positive number $a$.
 
The original dispersive estimates expressed $e^{itH}$ as a mapping between
weighted $L^2$ spaces, with the weights being exponential~\cite{Rauch} or
polynomial~\cite{JK}.  A significant advance was made by
Journ\'e, Soffer, and Sogge~\cite{JSS},
who proved the translation-invariant $L^1\to L^\infty$ bound 
in~\eqref{eq:dispersive} for potentials satisfying
$|V(x)| \le C(1+|x|)^{-7-\eps}$ and $\hat{V}\in L^1(\R^3)$.
The pointwise decay and regularity hypotheses were subsequently weakened 
by Yajima~\cite{Y1} and Goldberg and Schlag~\cite{GS1}, 
The ability to handle potentials with
$L^p$ singularities stems from recent results (e.g.~\cite{IonJer}) showing
that $-\Lap+V$ has no embedded eigenvalues at positive energies.

The hypotheses of Theorem~\ref{thm:dispersive} are nearly optimal
in a number of respects.  There exist compactly supported potentials
$V \in L^{3/2}_{\rm weak}$ for which $-\Lap + V$ admits bound states with
positive energy~\cite{KocTat}.  The inverse-square potential 
$V(x) = A|x|^{-2}$ only appears to be dispersive if
$A > -\frac14$~\cite{BPST}.
It is possible that the decay criteria can be relaxed slightly to include
all functions for which $\sup_y \int_{\R^3} |x-y|^{-1}|V(x)|\,dx$ is finite.
Such a condition is sufficient provided $V$ is small~\cite{RS}, or
mostly positive~\cite{DanPie}, and is critical with respect to scaling.
 
The proof of Theorem~\ref{thm:dispersive} begins by rewriting the
operator $e^{itH}P_{ac}(H)$ in terms of the resolvents 
$R_0(z) = (-\Lap-z)^{-1}$.  In this manner the dispersive estimate can be
reduced to a statement about the resolvents' mapping properties. 
As is frequently the case with dispersive phenomena, one needs to 
distinguish between high and low energies and make a separate calculation 
for each.  The various pieces are then assembled back into the original 
theorem at the end.
  
\subsection{Resolvent Identities}

Let $H=-\Lap+V$ in $\R^3$ and define the resolvents
$R_0(z) :=(-\Lap-z)^{-1}$ and $R_V(z):= (H-z)^{-1}$.
For $z \in \Compl \setminus \R^+$, the operator $R_0(z)$ can be realized as an
integral operator with the kernel
\[R_0(z)(x,y) = \frac{e^{i\sqrt{z}|x-y|}}{4\pi|x-y|} \]
where $\sqrt{z}$ is taken to have positive imaginary part.
While $R_V(z)$ does not possess an explicit representation of this form,
it can be expressed in terms of $R_0(z)$ via the identities
\begin{equation} \begin{aligned} \label{eq:res_ident}
R_V(z) &= (I + R_0(z)V)^{-1}R_0(z) = R_0(z)(I + VR_0(z))^{-1} \\
R_V(z) &= R_0(z) -R_0(z)VR_V(z)    = R_0(z) - R_V(z)VR_0(z)
\end{aligned}
\end{equation}
In the case where $z = \lambda \in \R^+$, one is led to consider limits
of the form $R_0(\lambda \pm i0) := \lim_{\eps\downarrow 0} 
R_0(\lambda\pm i\eps)$.  The choice of sign determines which branch of the
square-root function is selected in the formula above, therefore the two
continuations do not agree with one another.  
 For convenience we will adopt a shorthand
notation for dealing with resolvents along the positive real axis, namely
\begin{align*}
R_0^\pm(\lambda) &:= R_0(\lambda\pm i0) \\
R_V^\pm(\lambda) &:= R_V(\lambda\pm i0)
\end{align*}
Note that $R_0^-(\lambda)$ is the formal adjoint of $R_0^+(\lambda)$,
and a similar relationship holds for $R_V^\pm(\lambda^2)$.
The discrepancy between $R_0^+(\lambda)$ and $R_0^-(\lambda)$
characterizes the absolutely continuous part of the spectral measure
of $H$, denoted here by $E_{ac}(d\lambda)$, by means of the Stone formula
\begin{equation}
\label{eq:stone}
\langle E_{ac}(d\lambda) f,g \rangle
= \frac{1}{2\pi i}\big\langle [R_V^+(\lambda)-R_V^-(\lambda)]f,g \big\rangle 
\,d\lambda.
\end{equation}

Let $\chi$ be a smooth, even, cut-off function
on the line that is equal to one when $|x| \le 1$ and vanishes for all 
$|x| \ge 2$.  Further assume that translations of $\chi$ form a partition of
unity, in other words $\chi(x) + \chi(x-3) = 1$ for all $x \in [-1, 4]$. 
In order to prove Theorem~\ref{thm:dispersive}
it suffices to show that
\begin{align} 
\sup_{L\ge1}\Big|\Big \la &e^{itH} \chi^3(\sqrt{H}/L)\Pac f,g \Big\ra\Big|
\ =\   
\sup_{L\ge1}\Big|\int_0^\infty e^{it\lambda^2}\lambda\, \chi^3(\lambda/L) 
\Big\la [R_V^+(\lambda^2)-R_V^-(\lambda^2)]f,g \Big\ra
\, \frac{d\lambda }{\pi i}\Big| \nn \\
&=\ \sup_{L\ge1}\Big|\int_0^\infty e^{it\lambda^2}\lambda\, \chi^3(\lambda/L)
\Big\la \big[R_0^+(\lambda^2)(I + VR_0^+(\lambda^2))^{-1} - 
             R_0^-(\lambda^2)(I + VR_0^-(\lambda^2))^{-1}\big]f,g \Big\ra
\, \frac{d\lambda}{\pi i}\Big| \nn \\
&\les\  |t|^{-\frac32}\|f\|_1\|g\|_1. \label{eq:spec_theo}
\end{align}
The first equality is precisely \eqref{eq:stone}, and we have also made the
change of variable $\lambda \mapsto \lambda^2$. It is convenient to recall 
here that $R_0(z)$ is a holomorphic family of operators on the domain
$\Compl \setminus \R^+$, thus $R_0(z^2)$ is holomorphic on the upper 
half-plane.  Continuation onto the boundary $\{z = \lambda \in \R\}$ is
accomplished by taking limits from the interior.
\begin{equation*}
R_0^+(\lambda^2) := \lim_{\eps\to 0} R_0((\lambda+i\eps)^2)
  = \lim_{\eps\to 0} R_0(\lambda^2 + i\,\sign (\lambda)\eps)
\end{equation*}
For all $\lambda > 0$ this agrees with the previous definition of 
$R_0^+(\lambda^2)$, and for $\lambda < 0$ we have the identity
$R_0^+(\lambda^2) = R_0^-((-\lambda)^2)$.

Using this extended definition, the integral in \eqref{eq:spec_theo} can
be rewritten as
\begin{equation*} \tag{\ref{eq:spec_theo}'}
\sup_{L\ge1}\Big| \int_{-\infty}^\infty e^{it\lambda^2}\lambda\, 
\chi^3(\lambda/L) \big\la R_0^+(\lambda^2)(I + VR_0^+(\lambda^2))^{-1} f,g
\big\ra\, \frac{d\lambda}{\pi i} \Big|
\end{equation*}

{bf Remark.}
One can extend the domain of $R_0^-(\lambda^2)$ to the entire real line
by taking the domain of $R_0(z^2)$ to be the lower half-plane in $\Compl$.
The symmetry between $R_0^+(\lambda^2)$ and $R_0^-(\lambda^2)$ is reflected
in the identity
\begin{equation*}
R_0^-(\lambda^2) = R_0^+((-\lambda)^2) \quad {\rm for\ all}\ \lambda\in\R
\end{equation*}

It will be shown that, provided zero energy is not an eigenvalue or resonance,
the operators $T(\lambda) = (I + VR_0^+(\lambda^2))^{-1}$ are bounded on 
$L^1(\R^3)$, uniformly in $\lambda \in \R$.  We have included several copies
of the cutoff function $\chi(\lambda/L)$ in (\ref{eq:spec_theo}') so that
one of them may be combined with $T(\lambda)$ to form
\begin{equation} \label{eq:T_L}
T_L(\lambda) = \big(I + VR_0^+(\lambda^2)\big)^{-1} \chi(\lambda/L)
\end{equation}
For each $L\ge 1$ we have $\int_\R \norm[T_L(\lambda)][1\to 1]\, d\lambda
< \infty$,
therefore it has a well-defined Fourier transform with respect to $\lambda$.
Theorem~\ref{thm:dispersive} will eventually be derived from the following
estimate on the Fourier transform of $T_L$.

\begin{theorem} \label{thm:Fourier}
If $V$ satisfies the conditions of Theorem~\ref{thm:dispersive}, then
the family of operators $T_L(\lambda)$ have the property
\begin{equation}  \label{eq:Fourier}
\sup_{L\ge 1} \int_{\R^3}\int_\R |\widehat{T_L}(\rho)f(x)|\, d\rho\, dx
   \les \norm[f][1] \quad {\rm for\ all}\ f \in L^1(\R^3).
\end{equation}
\end{theorem}
\begin{prop}
The bound in \eqref{eq:Fourier} also holds for 
$T_L^-(\lambda) = (I + VR_0^-(\lambda^2))^{-1}\chi(\lambda/L)$.
\end{prop}
\begin{proof}
The family of operators $T_L^-(\lambda)$ is obtained from $T_L(\lambda)$
by taking complex conjugates, so $\widehat{T_L^-}(\rho)$ is just the complex
conjugate of $\widehat{T_L}(-\rho)$.  Neither conjugation nor reflection
changes the value of the inner integral over $\rho\in\R$.
\end{proof}

In the next two sections we will prove Theorem~\ref{thm:Fourier} by splitting
it into high-energy and low-energy cases.  For high energies, the argument
is a refinement of estimates found in \cite{RS} for each individual term
of the Born series.  The key step is a differentiability estimate which
enables us to control the geometric growth of the terms.  For low energies
the argument is an improvement of the one in \cite{GS1}, both in terms of
the computations required and the result achieved.  Finally, we show how
the dispersive bound in Theorem~\ref{thm:dispersive} follows from 
Theorem~\ref{thm:Fourier}.

\section{The High-Energy Case}

In this section we wish to show that Theorem~\ref{thm:Fourier} holds provided
we introduce a cutoff at sufficiently high energy.  A precise statement is
formulated below.
\begin{theorem} \label{thm:highenergy}
Let $V\in L^p(\R^3)\cap L^q(\R^3)$, $p<\frac32 < q$.
There exist a number $\lambda_1(V) < \infty$ and a constant
$A(V) < \infty$ so that the inequality
\begin{equation} \label{eq:high}
\sup_{L\ge 1} \int_{\R^3}\int_\R \big|
  \big[(1 - \chi(\cdot/\lambda_1))T_L\big]^\wedge(\rho)f(x)\big|\,d\rho\,dx 
     \le A \norm[f][1]
\end{equation}
holds for all $f \in L^1(\R^3)$.  
\end{theorem}
The general idea of the proof is to expand 
$T_L(\lambda) = \chi(\lambda/L)(I + VR_0^+(\lambda^2))^{-1}$ as a power series
and make estimates on each of the resulting terms.  The high-energy cutoff
will be needed only at the end to insure summability of the entire series.
We begin with an elementary observation.
\begin{proposition} \label{prop:Kato}
If $V \in L^p(\R^3)\cap L^q(\R^3)$, $p < \frac32 < q$, then
\begin{equation*}
\sup_{y \in \R^3} \int_{\R^3}\frac{|V(x)|}{|x-y|}\, dx \le C_{p,q}\norm[V][]
\end{equation*}
\end{proposition}
\begin{proof}
Here, and in the remainder of the discussion, we use $\norm[V][]$ to
indicate $\max(\norm[V][p], \norm[V][q])$.
Inside the region $\{|x-y| < 1\}$, use H\"older's inequality with 
$V \in L^q(\R^3)$ and $|\cdot - y|^{-1} \in L^{q'}(\R^3)$.  In the region
$\{|x-y| \ge 1\}$, consider $V \in L^p(\R^3)$ and $|\cdot - y|^{-1} \in 
L^{p'}(\R^3)$.
\end{proof}
\begin{corollary}
If $V \in L^p(\R^3) \cap L^q(\R^3)$, $p<\frac32<q$, then
\begin{equation*}
\norm[VR_0^+(\lambda^2)f][1] \le C_{p,q} \norm[f][1]
\quad {\it for\ all\ } f\in L^1(\R^3),\ \lambda \in \R
\end{equation*}
\end{corollary}
\begin{proof}
Recall that the free resolvent in three dimensions can be represented
explicitly by the integration kernel
\begin{equation} \label{eq:kernel}
R_0^+(\lambda^2)(x,y) = \frac{e^{i\lambda |x-y|}}{4\pi |x-y|}
\end{equation}
The integration kernel for $VR_0^+(\lambda^2)$ is therefore
$\frac{e^{i\lambda|x-y|}V(x)}{4\pi|x-y|}$, which immeditately satisfies the
Schur criterion for boundedness as an operator on $L^1$.
\end{proof}

\subsection{Integrability and Smoothness}
The next lemma is a fundamental $L^1$ estimate for the Fourier transform
of $VR_0^+(\lambda^2)$.  An unweighted version is implicit in the proof 
of Theorem~2.6 in \cite{RS}, however the extra decay assumption of $V$
(we have $V \in L^p(\R^3)$ instead of $\sup_y \int_{\R^3} |x-y|^{-1}|V(x)|\,dx
 < \infty$) allows
us to introduce a small polynomial weight. 
\begin{lemma}\label{lem:weightedL1}
 If $V\in L^p(\R^3)\cap L^q(\R^3)$, $p < \frac32 < q$, then there
exist $0< \eps < 1$ and $C < \infty$ such that
\begin{equation*}
\sup_{L\ge 1} \int_{\R^3}\int_\R \la\rho\ra^\eps \Big| \big[\chi(\cdot/L)
   \big(VR_0^+((\cdot)^2)\big)^k\big]^\wedge(\rho)f(x)\Big|\, d\rho\, dx
     \les (C\norm[V][])^k \norm[f][1]
\end{equation*}
for every $f \in L^1(\R^3)$ and every $k \ge 0$.
\end{lemma}
\begin{proof}
The expression $\la\rho\ra$ is given the usual meaning $(1 + |\rho|^2)^{1/2}$.
After substituting the integration kernel \eqref{eq:kernel} for each 
occurrence of $R_0^+(\lambda^2)$, we see that
\begin{equation*}
\chi(\lambda/L) \big(VR_0^+(\lambda^2)\big)^k f(x_0)
  = (4\pi)^{-k}  \int_{\R^{3k}} \chi(\lambda/L) e^{i\lambda\Sigma} 
 \left(\prod_{\ell = 0}^{k-1} \frac{V(x_{\ell})}{|x_{\ell} - x_{\ell+1}|}
 \right) f(x_k)\, dx_1dx_2 \ldots dx_k
\end{equation*}
where we have introduced the abbreviation $\Sigma := \sum_{\ell = 0}^{k-1}
 |x_\ell - x_{\ell+1}|$.  
The integrand above is a function in $L^1(\R^{3k+4})$.  This is most easily
seen by integrating sequentially in the variables $dx_0, dx_1,\ldots, dx_{k-1}$
and applying Proposition~\ref{prop:Kato} each time.
We may therefore use Fubini's theorem to take the Fourier
transform in $\lambda$ before integrating in $x_1, \ldots, x_k$.
The resulting expression is
\begin{equation*}
\big[\chi(\cdot/L)\big(VR_0^+((\cdot)^2)\big)^k\big]^\wedge (\rho) f(x_0)
 = (4\pi)^{-k}L \int_{\R^{3k}} \hat{\chi}(L(\rho - \Sigma))
 \left(\prod_{\ell = 0}^{k-1} \frac{V(x_{\ell})}{|x_{\ell} - x_{\ell+1}|} 
 \right) f(x_k)\, dx_1dx_2 \ldots dx_k
\end{equation*}

Multiply by the weight $\la\rho\ra^{\eps}$ and integrate with respect to
$d\rho$.  It is an elementary fact, proven below, that for any $\Sigma \in \R$,
$\sup_{L\ge 1}L \int_\R \la\rho\ra^\eps 
 \big|\hat{\chi}(L(\rho-\Sigma))\big|\, d\rho \les \la\Sigma\ra^\eps$.
Recall that $\Sigma = \sum_{\ell=0}^{k-1} |x_\ell - x_{\ell+1}|$ by definition.
Repeated application of the inequalities $\la A+B\ra < \la A\ra + \la B\ra$
and $(A+B)^\eps \le A^\eps + B^\eps$ for nonnegative $A,B$ and $0 < \eps < 1$
shows that $\la\Sigma\ra^\eps \le \sum_{\ell=0}^{k-1} 
\la|x_\ell - x_{\ell+1}|\ra^\eps$.
It follows that
\begin{multline*}
\sup_{L\ge 1}\int_R \la\rho\ra^\eps \big|\big[\chi(\cdot/L)
  \big(VR_0^+((\cdot)^2)\big)^k\big]^\wedge (\rho) f(x_0)\big|\,d\rho \\
\les (4\pi)^{-k}\sum_{\ell=0}^{k-1} \int_{\R^{3k}}
  \frac{|V(x_0)|}{|x_0-x_1|} \cdots 
  \frac{|V(x_\ell)|\la|x_\ell-x_{\ell+1}|\ra^\eps}{|x_\ell-x_{\ell+1}|} \cdots
  \frac{|V(x_{k-1})|}{|x_{k-1}-x_k|}|f(x_k)|\, dx_1\ldots dx_k
\end{multline*}
When the $L^1(R^3)$ norm is taken in the $x_0$ variable, it is possible to
integrate sequentially with respect to $dx_0, dx_1,\ldots dx_{k-1}$, applying
Proposition~\ref{prop:Kato} each time.  The integral in $dx_\ell$ is slightly
different, however it too is uniformly bounded provided 
$1-\eps > \frac3{p'}$.  Summing over $\ell$ introduces an extra factor of $k$,
however this may be absorbed into the constant $C$ because $k \le 2^k$ for
all $k \ge 0$.

By Fubini's theorem, the same result would be achieved had we first
integrated $dx_1\ldots dx_k$, then $d\rho$ and $dx$ as suggested in the
statement of the lemma.
\end{proof}
\begin{prop}
Let $\eta:\R\to\R$ satisfy the size bound $|\eta(y)| \les \la y\ra^{-2}$.
Then for any $0 < \eps < 1$,
\begin{equation*}
\sup_{L\ge 1}\ L \int_\R \la y\ra^\eps |\eta(L(y-\Sigma))| \les 
  \la\Sigma\ra^\eps
\end{equation*}
\end{prop}
\begin{proof}
If $|\Sigma| \le 1$, then the integral over the domain $y \in [-2,2]$ is 
comparable to 1, as desired.  For $|y| > 2$, $|\eta(L(y-\Sigma))| \les
|Ly|^{-2}$ Thus the tail of the integral is controlled by $\frac1L \le 1$.

If $|\Sigma| > 1$, integrate first on the domain $y \in [-2\Sigma, 2\Sigma]$,
taking $\la y\ra^\eps$ in $L^\infty$ and $L\eta(L(\cdot - \Sigma))$ in $L^1$.
Using the same estimates as above, the tail integral contributes no more
than $L^{-1}\Sigma^{\eps-1} \le \la\Sigma\ra^\eps$.
\end{proof}

The next order of business is to show that the Fourier transform of
$(VR_0^+(\lambda^2))^{k}$ becomes differentiable for sufficiently large $k$.
This corresponds to polynomial decay in $\lambda$ of $(VR_0^+(\lambda^2))^k$
as an operator on $L^1(\R^3)$.  We paraphrase the relevant statement from
\cite{Gol}.
\begin{proposition} \label{prop:lambdadecay} 
Let $V \in L^p(\R^3) \cap L^q(\R^3)$, $p < \frac32 < q$.  Then there exist
$\alpha > 0$ and $C < \infty$ such that
\begin{equation}
\norm[(VR_0^+(\lambda^2))^k f][1] \les (C\norm[V][])^k
   (1+\lambda)^{-(k-2)\alpha} \norm[f][1]
\end{equation}
\end{proposition}
\begin{proof}[Sketch of Proof]
It is a trivial matter to prove a uniform version of this bound, without any
deacy in $\lambda$, as the operator $VR_0^+(\lambda^2)$ is already known to
map $L^1(\R^3)$ to itself.  For $k = 0,1,2,$ this is sufficent.
The challenge is to use oscillation in the 
integration kernel $\frac{e^{i\lambda|x-y|}V(x)}{|x-y|}$ to strengthen the
bounds for large $\lambda$ and $k > 2$.

Choose a number $r \in \big(1, \min(\frac{3q}{q+3}, \frac{3p}{5p-3})\big)$.
The free resolvent $R_0^+(\lambda^2)$ is weak-type $(1,3)$, and also maps
$L^{\frac43}(\R^3)$ to $L^4(\R^3)$ with norm proportional to $\lambda^{-1/2}$
(This is a special case of Theorem~2.3 in \cite{KRS}).  By interpolation,
we conclude that $\norm[R_0^+(\lambda^2)][r\to 3r] = C_r\lambda^{-2/r'}$.
Because $\frac{3r}2 \in [p,q]$, it follows that
\begin{equation*}
\norm[VR_0^+(\lambda^2)][r\to r] \le C_r\norm[V][] \lambda^{-2/r'}
\end{equation*}
Two additional mapping estimates on $VR_0^+(\lambda^2)$ complete the proof.
Since $\bignorm[\frac{|V(\cdot)|}{|\cdot-y|}][r] < \infty$ uniformly in $y$, 
it is bounded as a map from $L^1(\R^3)$ to $L^r(\R^3)$.  By the 
Hardy-Littlewood-Sobolev theorem, it is also a bounded map from $L^r(\R^3)$
back to $L^1(\R^3)$.
\end{proof}
\begin{corollary} \label{cor:Sobolev}
Let $V\in L^p(\R^3) \cap L^q(\R^3)$, $p <\frac32 < q$.  then
\begin{equation*}
\sup_{L\ge 1}
 \int_{\R^3} \Big(\sup_{\rho\in\R} \big|(-\Lap_\rho+1)^{\frac{s}2}
  \big[\chi(\cdot/L)\big(VR_0^+((\cdot)^2)\big)^k\big]^\wedge(\rho)
  f(x)\big|\Big)\, dx 
   \les (C\norm[V][])^{k} \norm[f][1]
\end{equation*}
for any complex number $s$ with $\Re(s) < (k-2)\alpha - 1$.  The inequality
holds uniformly in $(s,k)$ satisfying $\Re(s) \le (k-2)\alpha - 2$.
\end{corollary}
\begin{proof}
By the previous lemma, $\la \lambda\ra^s (VR_0^+(\lambda^2))^k f$ will be
an integrable family of functions in $L^1(\R^3)$.  Multiplying by the cutoff
$\chi(\lambda/L)$ does not affect integrability.  By Fubini's theorem, 
that makes $(VR_0^+((\cdot)^2))^k f(x)$ an integrable
(over $x\in\R^3$) family of functions in $L^1(\la\lambda\ra^s\,d\lambda)$
Taken pointwise in~$x$, the Fourier transform in $\lambda$ maps this
space boundedly to $W^{s,\infty}$, as desired.
\end{proof}

\subsection{Interpolation}
We wish to interpolate between the estimates in Lemma~\ref{lem:weightedL1}
and Corollary~\ref{cor:Sobolev} to conclude that the Fourier transform
of $\chi(\cdot/L)(VR_0^+(\lambda^2))^kf$ has a small number of derivatives
in the space $L^1_xL^1_\rho$.  For technical reasons related to derivates of
imaginary order, it will be preferable to use $L^{1+\eps}$ with a polynomial
weight as a proxy for $L^1$.

As a preliminary step, observe that the case $s = 0$ in 
Corollary~\ref{cor:Sobolev} provides an $L^1_xL^\infty_\rho$ estimate on the 
function $\big[\chi(\cdot/L)(VR_0^+((\cdot)^2))\big]^\wedge(\rho)f(x)$,
while Lemma~\ref{lem:weightedL1} bounds its norm in $L^1_xL^1_\rho$.
Interpolate using H\"older's inequality to conclude that
\begin{equation*}
\sup_{L\ge 1}\int_{\R^3} \bignorm[\la\rho\ra^{\eps'}\big[\chi(\cdot/L)
 \big(VR_0^+((\cdot)^2)\big)^k\big]^\wedge(\rho)f(x)][L^{1+\eps/4}(d\rho)]\,dx
\les (C\norm[V][])^k \norm[f][1]
\end{equation*}
where $\eps' = \eps/(1+\frac{\eps}4)$.
The main step will be complex interpolation on the family of functions
\begin{equation*}
F_\theta(x,\rho) = \la\rho\ra^{\eps'(1-\theta)}(-\Lap_\rho+1)^{\frac{s}2\theta}
 \big[\chi(\cdot/L)\big(VR_0^+((\cdot)^2)\big)^k\big]^\wedge(\rho)f(x)
\end{equation*}
with $s = (k-2)\alpha - 2$ and
$\theta \in\Compl$ ranging over the strip $0 \le \Im(\theta) \le 1$.

On the boundary of the strip with $\theta = 1+i\gamma$, these functions are
uniformly bounded in $L^1_xL^\infty_\rho$ by Corollary~\ref{cor:Sobolev}
and the fact that $|\la\rho\ra^{-i\eps'\gamma}| = 1$.  For the boundary
with $\theta = i\gamma$, we use the fact that $(-\Lap_\rho+1)^{is\gamma/2}$
is a pseudodifferential operator of order zero, and can be represented by
convolution with a singular kernel $K_\gamma(\rho)$.  Following the 
calculations in \cite{bigStein}, chapter~6, one obtains the bounds
\begin{equation} \label{eq:C-Z}
|K_\gamma(\rho)| \les \la s\gamma\ra^2 |\rho|^{-1}, \qquad 
|K_\gamma'(\rho)| \les \la s\gamma\ra^3|\rho|^{-2}
\end{equation}
for all $x\in\R\setminus \{0\}$.  Additionally, since the second derivative
of $(1+\lambda^2)^{is\gamma/2}$ is integrable, $K_\gamma(x)$ satisfies the size
bound $|K_\gamma(x)| \les \la s\gamma\ra^{2}|x|^{-2}$.

For each value of $\gamma\in\R$, let $K_1(\rho) = \chi(\rho)K_\gamma(\rho)$
and $K_2(\rho) = (1-\chi(\rho))K_\gamma(\rho)$.  It is easy to verify that
$K_1(\rho)$ also satisfies the estimates in \eqref{eq:C-Z}, and that its
Fourier transform is a bounded function.  Then convolution with $K_1$ 
is a Calder\'on-Zygmund operator, hence it is bounded on $L^{1+\eps/4}(\R)$.
Moreover, since $K_1$ is supported on the interval $[-2,2]$,
\begin{equation*}
 \int_n^{n+1}\la\rho\ra^{\eps}|g*K_1(\rho)|^{1+\eps/4}\,d\rho 
   \les \la s\gamma\ra^{3(1+\eps/4)} \int_{n-2}^{n+3} \la\rho\ra^{\eps}
 |g(\rho)|^{1+\eps/4} \,d\rho
\end{equation*}
It is permissible to include the weight $\la\rho\ra^{\eps}$ is this inequality
because it has size comparable to $\la n\ra^{\eps}$ everywhere in both
domains of integration.  Summing over all $n\in\Z$,
\begin{equation*}
 \big(\norm[\la\rho\ra^{\eps'}g*K_1][1+\eps/4]\big)^{1+\eps/4} 
\les \la s\gamma\ra^{3(1+\eps/4)} \sum_{n\in\Z} \int_{n-2}^{n+3} 
 \la\rho\ra^{\eps}|g(\rho)|^{1+\eps/4}\,d\rho 
  = 5 \la s\gamma\ra^{3(1+\eps/4)} \big(\norm[\la\rho\ra^{\eps'}g][1+\eps/4]
   \big)^{1+\eps/4}
\end{equation*}

In other words, convolution with $K_1$ preserves the weighted space 
$L^{1+\eps/4}(\la\rho\ra^{\eps}d\rho)$.  The same is true of convolution
with $K_2$.  This is most readily seen by considering the action of
the integral kernel $\la\rho\ra^{\eps'}K_2(\rho-\sigma)\la\sigma\ra^{-\eps'}$
on unweighted $L^{1+\eps/4}(\R)$.  Note that
\begin{equation*}
\int_{|\rho-\sigma|>1} \frac{\la\rho\ra^{\eps'}}{|\rho-\sigma|^2}\,d\rho
  \les \la\sigma\ra^{\eps'}
\end{equation*}
If $|\sigma| < 2$ this is immediate. for $|\sigma| > 2$ break the domain
into the segments $\{|\rho| \le 2|\sigma|\}$ and $\{|\rho| > 2|\sigma|\}$.
Similarly, for any fixed $\rho\in\R$ we have
\begin{equation*}
\int_{|\rho-\sigma| > 1} \frac{\la\sigma\ra^{\eps'}}{|\rho-\sigma|^2}\,d\sigma
  \les \la\rho\ra^{-\eps'}
\end{equation*}
If $|\rho| < 2$ this is also immediate.  For $|\rho| > 2$, the domain of
integration should be broken into three pieces: $\{\sigma\in[\rho/2, 2\rho]\}$,
$\{\sigma\in[-2\rho,\rho/2]\}$, and $\{|\sigma| > 2|\rho|\}$.  Finally, one
concludes from the Schur test that convolution with $K_2$ is a bounded
operator on the weighted space $L^p(\la\rho\ra^{\eps'/p}d\rho)$ for any
exponent $1\le p\le\infty$, with particular emphasis on the case 
$p = 1+\eps/4$.  The operator norm is always less than $\la s\gamma\ra^2$,
regardless of the choice of $p$, since $|K_\gamma(x)| \les \la s\gamma\ra^2
|x|^{-2}$.

The end result of these calculations is that $(-\Lap_\rho+1)^{i\gamma}$ is 
bounded on the weighted space $L^{1+\eps/4}(\la\rho\ra^{\eps}d\rho)$, with 
operator norm growing at most polynomially in $|\gamma|$ and $s$.  
It follows that
\begin{equation*}
\norm[F_{i\gamma}][L^1_xL^{1+\eps/4}_\rho] \les \la s\gamma\ra^3
  (C\norm[V][1])^k \norm[f][1]
\quad {\rm and}\quad \norm[F_{1+i\gamma}][L^1_xL^\infty_\rho] \les 
 (C\norm[V][])^k \norm[f][1]
\end{equation*}
for all $\gamma\in\R$.  Apply complex interpolation and examine the case
$\theta = \frac{\eps}{4+2\eps}$.  The resulting bound is
\begin{equation*}
\bignorm[\la\rho\ra^{2\eps/(2+\eps)}(-\Lap_\rho+1)^{s\eps/(8+4\eps)}
 \big[\chi(\cdot/L)\big(VR_0^+((\cdot)^2)\big)^k\big]^\wedge(\rho)f(x)][L^1_x
  L^{1+\eps/2}_\rho] \les \la s\ra^3 (C\norm[V][])^k \norm[f][1]
\end{equation*}
The parameter $s$ was defined as a linear function of $k$, so the factor of
$\la s\ra^3$ may again be absorbed into the constant $C$ as in 
Lemma~\ref{lem:weightedL1}.
Observe that the reciprocal of $\la\rho\ra^{2\eps/(2+\eps)}$ is a function in 
$L^{(2+\eps)/\eps}(\R)$, the space dual to $L^{1+\eps/2}(\R)$.  H\"older's
inequality then leads to an estimate in $L^1_xL^1_\rho$, which we
formulate as a lemma.

\begin{lemma} \label{lem:interpolation}
Suppose $V \in L^p(\R^3)\cap L^q(\R^3)$, $p <\frac32 < q$, and let 
$k > \frac2\alpha + 2$.  Then
\begin{equation} \label{eq:interpolation}
\sup_{L\ge 1} \int_{\R^3} \int_\R \Big| (-\Lap+1)^{s_k}
 \big[\chi(\cdot/L)\big(VR_0^+((\cdot)^2)\big)^k\big]^\wedge(\rho)f(x) \Big|
 \, d\rho\,dx  \les (C\norm[V][])^k \norm[f][1]
\end{equation}
where $s_k = \frac{\alpha\eps}{8+4\eps}k - \frac{(1+\alpha)\eps}{4+2\eps}$
\end{lemma}

\subsection{Proof of Theorem~\ref{thm:highenergy}} 
Up to this point, our estimates have included the entire energy spectrum,
and the bounds grow geometrically in $k$
with ratio proportional to $\norm[V][]$.  The next lemma suggests how
introducing the high-energy cutoff $(1-\chi(\lambda/\lambda_1))$ can
ensure convergence of the geometric series even if $\norm[V][]$ is large.
\begin{lemma} 
Given $\lambda_1 > 1$ and $m > 0$, define a function
$F(\lambda)
= \la\lambda\ra^{-m}(1 - \chi(\lambda/\lambda_1))$.
\begin{equation*}
{\it Then}\ \norm[\hat{F}][1] \les \la m\ra\lambda_1^{-m}.
\end{equation*}
\end{lemma}
\begin{proof}
For $|\rho| > \lambda_1^{-1}$, use the identity 
\begin{equation*}
|\hat{F}(\rho)| = \Big| \rho^{-2} \Big(\frac{d^2F}{d\lambda^2}
\Big)^\wedge(\rho)\Big| 
 \le \rho^{-2} \bignorm[\frac{d^2F}{d\lambda^2}][1]
\end{equation*}
The second derivative can be computed using the product rule, and consists of
three terms.  Two of them are compactly supported on the intervals
where $|\lambda| \sim \lambda_1$, and are no larger than
$\la m\ra\lambda_1^{-(m+2)}$ anywhere on this set.  The last term, where both
derivates fall on $\la\lambda\ra^{-m}$, is supported where 
$|\lambda| \gtrsim \lambda_1$ and is everywhere smaller than
$\la m\ra^2\lambda^{-(m+2)}$.  The $L^1$ norm of each piece is seen to be less
than $(1+m) \lambda_1^{-(m+1)}$.  We conclude that
\begin{equation*}
|\hat{F}(\rho)| \les \la m\ra\lambda_1^{-(m+1)}\rho^{-2}
\end{equation*}
for all $|\rho| > \lambda_1^{-1}$.  This can contribute no more than
$\la m\ra \lambda_1^{-m}$ to the $L^1$-norm of $\hat{F}$.

If $m \ge \frac34$, then $\norm[F][2] \les \lambda_1^{-(m-1/2)}$.  By 
Plancherel's identity, the $L^2$-norm of $\hat{F}$ satisfies the same bound.
Then by the Cauchy-Schwartz inequality,
\begin{equation*}
\int_{-1/\lambda_1}^{1/\lambda_1} |\hat{F}(\rho)|\, d\rho \les 
\lambda_1^{-m}
\end{equation*}

For $m < \frac34$, write $\la\lambda\ra^{-m}(1-\chi(\lambda/\lambda_1)) = 
|\lambda|^{-m} + G(\lambda)$.  The remainder function $G$ is dominated by
$|\lambda|^{-m}$ for all $\lambda \le 2\lambda_1$ and by $m|\lambda|^{-(m+2)}$
for all $\lambda > 2\lambda_1$.  Thus 
$\norm[G][1] \les \la m\ra \lambda_1^{1-m}$.

The Fourier transform of $|\lambda|^{-m}$ is exactly $c_m|\rho|^{m-1}$, where
$c_m \les m$ for the range of $m$ under consideration.  The $L^1$-norm of
this function on the interval $[-\lambda_1^{-1}, \lambda_1^{-1}]$ is 
$\frac{c_m}{m} \lambda_1^{-m}$, and the constant $\frac{c_m}{m}$ is bounded
uniformly.  Meanwhile, the Fourier transform of $G$ is bounded above by
$\la m\ra \lambda_1^{1-m}$, so its $L^1$-norm over the same interval is
less than $\la m\ra \lambda_1^{-m}$ as well.  Adding the two pieces together
proves the desired estimate. 
\end{proof}

\begin{proof}[Proof of Theorem~\ref{thm:highenergy}]
Recall that we are trying to verify the inequality 
\begin{equation*} 
\sup_{L\ge 1} \int_\R \bignorm[\big[(1 - \chi(\cdot/\lambda_1))T_L\big]^\wedge
  (\rho)f][1]\, d\rho  \les \norm[f][1]     \tag{\ref{eq:high}}
\end{equation*}
Fix $L \ge 1$, and assume that $\lambda_1^\alpha > 2C\norm[V][]$, where 
$\alpha$ and $C$ are the
constants in Proposition~\ref{prop:lambdadecay}.  The power series
\begin{equation} \label{eq:TLseries}
(1 - \chi(\lambda/\lambda_1))T_L(\lambda) = \sum_{k=0}^\infty
  (1 - \chi(\lambda/\lambda_1)) \chi(\lambda/L)\big(VR_0^+(\lambda^2)\big)^k
\end{equation}
converges uniformly in $\lambda$, and is supported on the interval
$\lambda \in [\lambda_1, L]$.  The Fourier transform of the partial sums
then converges in the sense of distributions.

For $k \le \frac2\alpha + 2$, we use the estimate in
Lemma~\ref{lem:weightedL1} showing that
$\big[\chi(\cdot/L)\big(VR_0^+((\cdot)^2)\big)^k\big]^\wedge(\rho) f$ is
an integrable family (indexed by $x\in\R^3$) of functions in $L^1(\R)$.  
The Fourier transform of $(1 - \chi(\cdot/\lambda_1))$ is a measure whose 
total variation norm is finite and does not depend on $\lambda_1$.   
Each of these terms then contributes no more than $(C\norm[V][])^k\norm[f][1]$
to the total on the right-hand side of \eqref{eq:high}.

For all $k > \frac2\alpha + 2$, multiply and divide the $k^{\rm th}$ term 
by a common factor to obtain
\begin{equation*}
\Big(\la\lambda\ra^{-2s_k}(1-\chi(\lambda/\lambda_1))\Big)
   \Big(\la\lambda\ra^{2s_k} \chi(\lambda/L)\big(VR_0^+(\lambda^2)\big)^k\Big)
\end{equation*}
where $s_k$ is the same number as in~\eqref{eq:interpolation}.  
Consider the second factor in this product.  
By Lemma~\ref{lem:interpolation} its Fourier transform, acting on a fixed
function $f$, also gives rise to an integrable family of $L^1(\R)$ functions
indexed by $x\in\R^3$.  The $L^1$-norm of this family is bounded by 
$(C\norm[V][])^k \norm[f][1]$.

The Fourier transform of the first factor is an integrable function of $\rho$,
with norm less than $\la s_k\ra \lambda_1^{-2s_k}$.  When this is convolved
against the expression from the second factor, the result is again an 
integrable family of $L^1(\R)$ functions with the norm bound
\begin{equation}
\int_{\R^3} \int_\R \Big|\big[(1-\chi(\cdot/\lambda_1))\chi(\cdot/L)
   \big(VR_0^+((\cdot)^2)\big)^k\big]^\wedge(\rho)f(x)\Big|\,d\rho\,dx
  \les \la s_k\ra (C\norm[V][])^k \lambda_1^{-2s_k} \norm[f][1]
\end{equation}
holding uniformly in $L$.  Recall that 
$s_k = \frac{\alpha\eps}{8+4\eps}k - \frac{(1+\alpha)\eps}{4+2\eps}$ by
definition, so $s_k$ is a linear function of $k$.  The bound shown above
is then geometric in $k$, and its ratio is moderated by a negative power
of $\lambda_1$.
If $\lambda_1$ is chosen so that $\lambda_1^{(\alpha\eps)/(4+2\eps)} > 
 2C\norm[V][]$, the geometric series converges, therefore
\begin{equation*}
\sup_{L\ge1} \sum_{k=0}^\infty \int_{\R^3}\int_\R
\big|\big[(1-\chi(\cdot/\lambda_1))\chi(\cdot/L)
   \big(VR_0^+((\cdot)^2)\big)^k\big]^\wedge(\rho) f(x)\big|\,d\rho\,dx
  \les (C\norm[V][])^{2+2/\alpha}\norm[f][1]
\end{equation*}
by comparing the entire series to its largest term.

The Fourier transform of the series \eqref{eq:TLseries} converges
in $L^1$ as well as in the distributional sense, and its limit has norm
controlled by $A(V)\norm[f][1]$.
\end{proof}

\section{The Low-Energy Case} \label{sec:lowenergy}
In this section we prove the complementary statement to 
Theorem~\ref{thm:highenergy}, namely
\begin{theorem} \label{thm:lowenergy}
Let $V$ satisfy the conditions of Theorem~\ref{thm:dispersive}, and fix any
$0 < \lambda_1 < \infty$
\begin{equation} \label{eq:low} 
\sup_{L\ge 1} \int_{\R^3}\int_\R \big|\big[\chi(\cdot/\lambda_1)T_L\big]^\wedge
  (\rho)f(x)\big|\,d\rho\,dx  \les \norm[f][1]
\end{equation} 
holds for all $f \in L^1(\R^3)$. 
\end{theorem}
There are two low-energy cutoffs present in the statement of this theorem,
since $T_L$ is shorthand for $\chi(\lambda/L)(I + VR_0^+(\lambda^2))^{-1}$.
We will relegate both of these to the background by introducing a third
cutoff function which localizes to much smaller intervals in $\lambda$.
The theorem is then proved by adding up a finite number of local results.

At low energies the Neumann series expansion of $(I + VR_0^+(\lambda^2))^{-1}$
will typically diverge unless $\norm[V][]$ is small.  The existence of 
inverses must instead be demonstrated by a Fredholm alternative argument.
For $\lambda = 0$ this requires that zero energy is neither an eigenvalue
nor a resonance, as defined below.
\begin{definition}
We say that a resonance occurs at zero energy if the equation 
$(I+VR_0^+(0))g = 0$ admits a distributional solution $g \not\in L^2(\R^3)$
such that $\la x\ra^{-\beta}g \in L^2(\R^3)$ for every $\beta > \frac12$.
\end{definition}

The Fredholm alternative does not construct inverses explicitly, 
which limits our ability to perform subsequent calculations.
We therefore avoid its use, except in a finite number of instances, by the 
following scheme:

Fix a ``benchmark'' energy $\lambda_0 \in \R$ and let 
$S_0 = (I+VR_0^+(\lambda_0^2))^{-1}$.  For all values of $\lambda$ 
sufficiently close to $\lambda_0$, we may regard $R_0^+(\lambda^2)$ as 
a perturbation of $R_0^+(\lambda_0^2)$ and treat the corresponding inverse
as a perturbation of $S_0$.  The underlying perturbation is a difference of
free resolvents, hence it can be represented explicitly by an integration 
kernel.  The role of $S_0$ is limited to its existence as a (fixed) bounded
operator on $L^1(\R^3)$.
In this manner the entire interval of energies $|\lambda -\lambda_0| < \delta$
may be considered with only one application of the Fredholm theory.
The perturbation radius $\delta$ can be chosen independent of $\lambda_0$,
so the low-energy spectrum $\lambda \in [-2\lambda_1, 2\lambda_1]$ is
covered by a finite collection of such intervals.

The details of the proof are clearly foreshadowed by the low-energy
discussion in \cite{GS1}.  Two technical modifications allow us to work with
a larger class of potentials while reducing the burden of computation.
One is the use of $L^1(\R^3)$ as the natural setting instead of weighted
$L^2$ spaces.  The other is, for a family of operators $T(\rho)$, 
estimating the quantity $\int_\R \norm[T(\rho)f][]\,d\rho$ rather than
$\int_\R \norm[T(\rho)][]\,d\rho$.  

\subsection{Invertibility of $I + VR_0^+(\lambda^2)$}
Here we show that $(I + VR_0^+(\lambda^2))^{-1}$ exists as a bounded operator
on $L^1(\R^3)$ for each $\lambda\in\R$, and that the operator norm of these
inverses can be controlled uniformly in $\lambda$.  Essentially identical
arguments have been made in various function spaces, and with varying 
assumptions on $V$, for example in~\cite{DanPie} and~\cite{GS2}, and can
traced back to Agmon's work on the limiting absorption principle~\cite{Agmon}.

\begin{lemma}
Suppose $V$ satisfies the conditions of Theorem~\ref{thm:dispersive}.
Then
\begin{equation} \label{eq:supinverse}
\sup_{\lambda\in\R} \norm[(I + VR_0^+(\lambda^2))^{-1}][1\to 1] < \infty
\end{equation}
\end{lemma}
\begin{proof}[Sketch of Proof]
Observe that if $V \in C^\infty_c(\R^3)$, then $VR_0^+(\lambda^2)$ maps
$L^1(\R^3)$ to $W^{2,1}(\supp V)$, hence it is a compact operator on 
$L^1(\R^3)$ by Rellich's theorem.  For general potentials, compactness of
$VR_0^+(\lambda^2)$ is seen by writing $V$ as a limit of functions in
$C^\infty_c(\R^3)$.

The Fredholm Alternative Theorem then dictates that either 
$(I + VR_0^+(\lambda^2))^{-1}$ is bounded on $L^1(\R^3)$ or else it has 
a nonempty null-space, that is there exists $g_\lambda \in L^1(\R^3)$ solving
$(I + VR_0^+(\lambda^2))g_\lambda = 0$.  By bootstrapping the identity
$g_\lambda = -VR_0^+(\lambda^2)g_\lambda$ with the 
Hardy-Littlewood-Sobolev inequality,
we see that $g_\lambda \in L^1(\R^3)\cap L^q(\R^3)$ and 
$R_0^+(\lambda^2)g_\lambda \in L^\infty(\R^3)$.  Since
\begin{equation*}
0 = \Im \la R_0^+(\lambda^2)g_\lambda, VR_0^+(\lambda^2)g_\lambda\ra = 
   -\Im \la R_0^+(\lambda^2)g_\lambda, g_\lambda\ra =  
   c\lambda\int_{S^2} |\hat{g_\lambda}(\lambda \omega)|^2\,d\omega
\end{equation*}
it follows that for any $\lambda \in \R \setminus \{0\}$, the Fourier 
transform of $g_\lambda$ vanishes (in the $L^2$ trace sense) on the sphere 
of radius $\lambda$.  By Proposition~12 in~\cite{GS2}, 
$R_0^+(\lambda^2)g_\lambda \in L^2(\R^3)$.  On the other hand, the definition
of $g_\lambda$ implies that
$(-\Lap + V - \lambda^2)(R_0^+(\lambda^2)g_\lambda) = 0$ in the sense of 
distributions.  Finally, a theorem of Ionescu and Jerison~\cite{IonJer} states
that $(-\Lap + V)$ has no nontrivial eignefunctions with positive energy,
so $R_0^+(\lambda^2)g_\lambda = 0$.  It follows immediately that 
$g_\lambda = 0$ as well.

In the case $\lambda = 0$, the distributional equation 
$(-\Lap + V)R_0^+(0)g_0 = 0$ is still valid.
For any $\beta > \frac12$, $\sup_{y\in\R^3} 
  \norm[\la x\ra^{-2\beta}|x-y|^{-1}][2] < \infty$, therefore
$R_0^+(0)$ is a bounded map from $L^1(\R^3)$ to the weighted space
$L^{2,-\beta}(\R^3)$.  By our assumption that zero energy is neither
an eigenvalue nor a resonance, we exclude the possiblity that a nontrivial
function $R_0^+(0)g_0$ can belong to this class, leaving only the
solution $g_0 = 0$.

So far we have established that $(I + VR_0^+(\lambda^2))^{-1}$ exists
at each $\lambda \in \R$, but have not shown uniformity.  By 
Proposition~\ref{prop:lambdadecay}, $\norm[(VR_0^+(\lambda^2))^3][1\to 1]
 < \frac12$ for sufficiently large $\lambda$.  For these values of $\lambda$,
\begin{equation*}
\norm[(I+ VR_0^+(\lambda^2))^{-1}][1\to1] \le \norm[I - VR_0^+(\lambda^2)
  + (VR_0^+(\lambda^2))^2][1\to1] \
   \norm[(I + (VR_0^+(\lambda^2))^3)^{-1}][1\to1]
  \les 1 + \norm[V][]^2
\end{equation*}
which provides a uniform bound.  For small $\lambda$, observe that the family
of operators $I + VR_0^+(\lambda^2)$ vary continuously in $\lambda$
(In fact, the variation is H\"older continuous because $V \in L^p(\R^3)$ for 
some $p < \frac32$).  Since inverses exist at every $\lambda \in \R$, they
also form a continuous family of operators, and are uniformly bounded on
any compact set.
\end{proof}

\subsection{Proof of Theorem~\ref{thm:lowenergy}}
Fix $\lambda_0 \in \R$ and let $S_0 = (I + VR_0^+(\lambda_0^2))^{-1}$.
Then $I + VR_0^+(\lambda^2) = S_0^{-1} + VB^+(\lambda)$, where $B^+(\lambda)$
denotes the  difference $R_0^+(\lambda^2) - R_0^+(\lambda_0^2)$.
Taking inverses,
\begin{equation} \label{eq:SVB}
\big(I + R_0^+(\lambda^2)\big)^{-1} = \big(I + S_0VB^+(\lambda)\big)^{-1}S_0
\end{equation}

We remarked above that $VR_0^+(\lambda^2)$ varies continuously in $\lambda$,
which suggests that $VB^+(\lambda)$ should vanish in the limit 
$\lambda \to \lambda_0$.  More precisely, $VB^+(\lambda)$ is an integral
operator with associated kernel
\begin{equation*}
|VB^+(\lambda,x,y)| = \left|\frac{V(x)\big(e^{i\lambda|x-y|}-
  e^{i\lambda_0|x-y|}\big)}{|x-y|}\right|  \le \left\{
\begin{aligned}
|\lambda&-\lambda_0|\,|V(x)|, & &{\rm if}\ |x-y| \le |\lambda-\lambda_0|\\ 
&\frac{|V(x)|}{|x-y|}, & &{\rm if}\ |x-y| > |\lambda-\lambda_0|
\end{aligned}  \right.
\end{equation*}
For fixed $y\in \R^3$, $\lambda\in\R$, the $L^1$-norm of this kernel in
the $x$ variable is controlled by $|\lambda-\lambda_0|^{1-3/p'}\norm[V][]$,
by applying H\"older's inequality with $V\in L^{p}(\R^3)$ and the
remaining factors in $L^{p'}(\R^3)$.  In other words, 
$\norm[VB^+(\lambda)][1\to1] \le C|\lambda-\lambda_0|^{1-3/p'}\norm[V][]$,
with the constant $C < \infty$ independent of the choice of $\lambda_0$.
Since $\norm[S_0][1\to1]$ is bounded above by \eqref{eq:supinverse},
there exists $r>0$ so that $\norm[S_0VB^+(\lambda)][1\to1] < \frac12$
whenever $|\lambda-\lambda_0| \le 4r$.

In that case, the Neumann series
\begin{equation} \label{eq:Neumann}
\chi(\frac{\lambda-\lambda_0}{r})\big(I+VR_0^+(\lambda^2)\big)^{-1} =
\sum_{k=0}^\infty (-1)^k\Big(\chi(\frac{\lambda-\lambda_0}{2r})
  \big(S_0VB^+(\lambda)\big)\Big)^k \chi(\frac{\lambda-\lambda_0}{r})S_0
\end{equation}
converges uniformly over all $\lambda\in\R$.  Recall here that 
$\chi(\frac{\lambda-\lambda_0}{2r}) = 1$ everywhere on the support of 
$\chi(\frac{\lambda-\lambda_0}{r})$.  The Fourier transforms of
the partial sums converge in the sense of distributions.

\begin{lemma}
The Fourier transform of $\chi(\frac{\lambda-\lambda_0}{2r})S_0VB^+(\lambda)$
satisfies the bound
\begin{equation}
\int_\R \bignorm[\big[\chi(\frac{\cdot-\lambda_0}{2r})S_0VB^+(\cdot)
  \big]^\wedge(\rho) f][1]\, d\rho \le C r^{1 -3/p'}\norm[f][1]
\end{equation}
for all functions $f \in L^1(\R^3)$.
\end{lemma}
\begin{proof}
The Fourier transform of $\chi(\frac{\lambda-\lambda_0}{2r})VB^+(\lambda)$ can
be represented by the integration kernel
\begin{equation*}
K_r(\rho,x,y) = e^{-i\lambda_0(\rho-|x-y|)}V(x)
   \left(2r\frac{\hat{\chi}(2r(\rho-|x-y|))-\hat{\chi}(2r\rho)}{|x-y|}\right)
\end{equation*}
This leads to an immediate estimate $\int_\R |K_r(\rho,x,y)|\,d\rho
 \le 2|x-y|^{-1}|V(x)|\norm[\hat{\chi}][1]$, by assuming no cancellation
 between the two evaluations of $\hat{\chi}$.  For small values of $|x-y|$
 a better estimate is possible.  By the Mean Value theorem,
\begin{equation*}
|K_r(\rho,x,y)| \le |V(x)| \frac{2r}{|x-y|}\int_{2r(\rho-|x-y|)}^{2r\rho}
  |\hat{\chi}'(\tau)|\,d\tau
\end{equation*}
Fubini's theorem permits integrations to be carried out in any order, so that
$\int_\R |K_r(\rho,x,y)|\,d\rho \le 2r|V(x)|\,\norm[\hat{\chi}'][1]$.
Putting the two estimates together,
\begin{equation*}
\int_\R |K_r(\rho,x,y)|\,d\rho \les |V(x)|\min(r,|x-y|^{-1})
\end{equation*}
which leads to the further integral estimate
\begin{equation*}
\iiint_{\R^{1+3+3}} |K_r(\rho,x,y)|\,f(y)\, d\rho\,dx\,dy
  \les r^{1-3/p'}\norm[V][] \norm[f][1].
\end{equation*}
Once again, Fubini's theorem allows for the integration to take place
in the reverse order.  This means that $\int_\R 
\norm[[\chi((\cdot-\lambda_0)/2r)VB^+(\cdot)]^\wedge(\rho)f][1]\,d\rho
  \les r^{1-3/p'}\norm[f][1]$.  Applying the bounded operator $S_0$ pointwise
at each $\rho \in \R$ only increases the estimate by a finite factor.
\end{proof}

\begin{corollary}
If $f_\rho$ is a family of functions in $L^1(\R^3)$ indexed by $\rho\in\R$,
then
\begin{equation}
\int_\R \bignorm[\int_\R \big[\chi(\frac{\cdot-\lambda_0}{2r})S_0B^+(\cdot)
  \big]^\wedge(\sigma - \rho) f_\rho\, d\rho][L^1(\R^3)] \,d\sigma
\le Cr^{1-3/p'} \int_\R \norm[f_\rho][1]\,d\rho
\end{equation}
\end{corollary}
\begin{proof}
The expression on the left-hand side is dominated by
\begin{equation*}
\iint_{\R^2} \bignorm[ \big[\chi(\frac{\cdot-\lambda_0}{2r})S_0VB^+(\cdot)
  \big]^\wedge(\sigma-\rho) f_\rho][1]\, d\rho\, d\sigma
\end{equation*}
which, after applying Fubini's theorem and the previous lemma, is seen to be
less than the expression on the right-hand side.
\end{proof}
A pointwise product of functions in the $\lambda$ variable corresponds to 
convolution in the $\rho$ variable when Fourier transforms are taken.
The previous two statements can be combined iteratively to prove Fourier
bounds for $(S_0VB^+(\lambda))^k$.
\begin{corollary}
The Fourier transform of $\big(\chi(\frac{\lambda-\lambda_0}{2r})S_0
  VB^+(\lambda)\big)^k$ satisfies the bound
\begin{equation} 
\int_\R \bignorm[\big[\big(\chi(\frac{\cdot-\lambda_0}{2r})S_0VB^+(\cdot)
 \big)^k \big]^\wedge(\rho) f][1]\, d\rho \le (Cr^{(1-3/p')})^k \norm[f][1]
\end{equation}
for all functions $f \in L^1(\R^3)$.
\end{corollary}
\begin{proof}[Proof of Theorem~\ref{thm:lowenergy}]
Apply the corollary above to $S_0f \in L^1(\R^3)$, then convolve 
in $\rho$ with the function $r\hat{\chi}(r\rho)$.
This has $L^1$-norm $\norm[\hat{\chi}][1] < \infty$, and $S_0$ is a bounded 
map, so the previous estimates are multiplied by a fixed constant.
If $r>0$ is chosen small enough so that $Cr^{1-3/p'} < \frac12$, then the
Neumann series for $\chi(\frac{\cdot-\lambda_0}{r})(I+VR_0^+(\lambda^2))^{-1}$
given in \eqref{eq:Neumann}
 converges in $L^1$-norm (as well as in distributions) on the Fourier 
transform side, and has norm bounded by $\norm[f][1]$.  Note that the 
chosen value of $r$ does not depend on $\lambda_0$.

Further convolutions in $\rho$ with the functions 
$\lambda_1\hat{\chi}(\lambda_1\rho)$ and $L\hat{\chi}(L\rho)$, 
each of which also has $L^1$-norm $\norm[\hat{\chi}][1]$, yields a similar 
estimate for $\big[\chi(\frac{\lambda-\lambda_0}{r})\chi(\lambda/\lambda_1)
  T_L\big]^\wedge$.  Since translations of $\chi$ form a partition of unity,
the localization caused by $\chi(\frac{\lambda-\lambda_0}{r})$ may be
removed by obtaining separate bounds for each choice of $\lambda_0 = 3nr$,
$n \in [-\frac{2\lambda_1}r, \frac{2\lambda_1}r]$, and adding these together.
\end{proof}

\section{Proof of Theorem~\ref{thm:dispersive}}
We now return to the goal of proving
\begin{equation*}
\sup_{L\ge1}\Big|\int_0^\infty e^{it\lambda^2}\lambda\, \chi^3(\lambda/L)
\Big\la [R_V^+(\lambda^2)-R_V^-(\lambda^2)]f,g \Big\ra
\, \frac{d\lambda }{\pi i}\Big|
\les |t|^{-\frac32} \norm[f][1] \norm[g][1] \tag{\ref{eq:spec_theo}}
\end{equation*}
Integrate the left-hand expression by parts once to obtain
\begin{equation} \label{eq:IBP}
\begin{aligned}
\frac1{2\pi|t|} \sup_{L\ge1} \Big|\int_0^\infty e^{it\lambda^2}
&\Big(\chi^3(\lambda/L)\Big\la \big[\dfrac{d}{d\lambda}R_V^+(\lambda^2) -
  \dfrac{d}{d\lambda}R_V^-(\lambda^2)\big]f, g\Big\ra \\
  &+\ \frac3L\, \chi'(\lambda/L)\chi^2(\lambda/L)\big\la[R_V^+(\lambda^2)
   - R_V^-(\lambda^2)]f,g\big\ra \Big)\,d\lambda \Big|
\end{aligned}
\end{equation}
The two terms are considered separately.  For the first one, 
use the identity $R_V^+((-\lambda)^2) = R_V^-(\lambda^2)$ to rewrite it as
\begin{equation*}
\frac1{2\pi|t|} \sup_{L\ge1} \Big|\int_{-\infty}^\infty e^{it\lambda^2}
 \chi^3(\lambda/L)\Big\la \dfrac{d}{d\lambda}R_V^+(\lambda^2)f,g\Big\ra\,
 d\lambda
\end{equation*}
compare this to (\ref{eq:spec_theo}').  The operator $\dfrac{d}{d\lambda}
R_V^+(\lambda^2)$ can be written in terms of free resolvents by 
differentiating the identity \eqref{eq:res_ident}.  There are several
algebraically equivalent expressions to choose from, one of which is
$(I + R_0^+(\lambda^2)V)^{-1}\dfrac{d}{d\lambda}\big[R_0^+(\lambda^2)\big]
(I+ VR_0^+(\lambda^2))^{-1}$.
Substituting this into the integral yields
\begin{equation*}
\frac1{2\pi|t|} \sup_{L\ge1} \Big|\int_\R e^{it\lambda^2}
\Big\la \chi(\lambda/L)\dfrac{d}{d\lambda}\big[R_0^+(\lambda^2)\big]
  T_L(\lambda)f, T_L^-(\lambda)g \Big\ra\, d\lambda
\end{equation*}
We wish to apply Parseval's theorem, separating the integrand into the
product $e^{it\lambda^2} \cdot A(\lambda)$.  The factor denoted by 
$A(\lambda)$ is bounded with compact support, since every operator
$T_L(\lambda)$ and $T_L^-(\lambda)$ is bounded on $L^1(\R^3)$ and
$\dfrac{d}{d\lambda}R_0^+(\lambda^2)$ maps $L^1(\R^3)$ to $L^\infty(\R^3)$.
More precisely, it has integral kernel $(-4\pi i)^{-1}e^{i\lambda|x-y|}$.
Thus the Fourier transform of $\chi(\lambda/L)\dfrac{d}{d\lambda}
 R_0^+(\lambda^2)$ is a family of integral operators with kernel
\begin{equation*}
   K(\rho,x,y) = \frac{iL}{2}\hat{\chi}(L(\rho-|x-y|)).
\end{equation*}
The Fourier transform of $e^{it\lambda^2}$ is well known to be
$\sqrt{\pi/(2|t|)}(1 + i\;\sign(t)) e^{(-i\rho^2/4t)}$.  Thus Parseval's 
theorem leads us to evaluate
\begin{equation}
|t|^{-\frac32} \sup_{L\ge1} \Big| L \int_\R e^{(-i\rho^2/4t)} 
\iiiint_{\R^8}  \big[\widehat{T_L}(\sigma)f(x)\big] 
  \big[\overline{\widehat{T_L^-}(\tau)g(y)}\big] 
  \hat{\chi}(L(\rho-\sigma-\tau-|x-y|))\,dx\,d\sigma\,dy\,d\tau\ d\rho
\end{equation}
modulo constants.  The fact that $A(\lambda)$ is a product of three terms
means that $\hat{A}(\rho)$ is an iterated convolution, hence the presence of
auxilliary variables $\sigma$ and $\tau$.  Take the absolute value inside
all the integrals, so that we may evaluate them in a more convenient order.

The integral $\int_\R L |\hat{\chi}(L(\rho-\sigma-\tau-|x-y|))|\,d\rho$
contributes $\norm[\hat{\chi}][1]$ for any fixed value of the other
variables.  Then, since $\widehat{T_L}(\sigma)f$ and
$\widehat{T_L^-}(\tau)g$ are both integrable families of functions in 
$L^1(\R^3)$ by Theorem~\ref{thm:Fourier},
the entire expression is controlled by $|t|^{-3/2} \norm[f][1]\norm[g][1]$.

The second term in \eqref{eq:IBP} is treated similarly.  $R_V^+(\lambda^2)
- R_V^-(\lambda^2)$ becomes an odd function when both pieces are extended
to all of $\lambda\in\R$, as is $\chi'(\lambda/L)$, so the entire integrand
is even.  We can then evaluate
\begin{equation*}
\frac3{4\pi|t|} \sup_{L\ge 1} \frac1L \Big|\int_{-\infty}^\infty 
  \chi'(\lambda/L)\chi^2(\lambda/L)\big\la[R_V^+(\lambda^2)
   - R_V^-(\lambda^2)]f,g\big\ra\, d\lambda \Big|
\end{equation*}
Cancellation between the two resolvents plays a much greater role here,
as we need $R_V^+(\lambda^2) - R_V^-(\lambda^2)$ to map $L^1(\R^3)$ into
$L^\infty(\R^3)$ for the inner product to be well-defined.  Starting with
the relations $\R_V^\pm(\lambda^2) = R_0^\pm(\lambda^2)
 (I + VR_0^\pm(\lambda^2))^{-1}$ and performing some algebra, we obtain
the identity
\begin{equation*}
R_V^+(\lambda^2) - R_V^-(\lambda^2) = \big(I+R_0^-(\lambda^2)V\big)^{-1}
   \big(R_0^+(\lambda^2) - R_0^-(\lambda^2)\big)
   \big(I+VR_0^+(\lambda^2)\big)^{-1}
\end{equation*}
The middle factor is precisely convolution with the kernel
$\frac{-\sin(\lambda|x|)}{2\pi |x|}$, which indeed maps $L^1(\R^3)$ to 
$L^\infty(\R^3)$, and the outer factors are each bounded on their respective
spaces.  Written another way, the expression in question is
\begin{equation*}
\frac{3}{4\pi|t|}\sup_{L\ge 1}\frac1L \Big|\int_\R
\big\la \chi'(\lambda/L)\big(R_0^+(\lambda^2)-R_0^-(\lambda^2)\big)
  T_L(\lambda)f, T_L(\lambda)g\Big\ra\, d\lambda \Big|
\end{equation*}
After applying Plancherel's theorem and discarding fixed constants, this is
equivalent to
\begin{multline} \label{eq:boundary}
|t|^{-\frac32}\sup_{L\ge 1} \Big|\int_\R e^{(-i\rho^2/4t)}
  \iiiint_{\R^8} \big[\widehat{T_L}(\sigma)f(x)\big]
  \big[\overline{\widehat{T_L}(\tau)g(y)}\big] \\
\times \frac{\widehat{\chi'}(L(\rho-\sigma-\tau -|x-y|)) - 
   \widehat{\chi'}(L(\rho-\sigma-\tau+|x-y|))}{|x-y|}
   \,dx\,d\sigma\,dy\,d\tau\ d\rho
\end{multline} 
Take the absolute value inside all integrals, and perform integration
first with respect to $d\rho$.  The complex exponential function disppears,
and we are left with
\begin{equation*}
\frac1{|x-y|} \int_\R \Big| \int_{L(\rho-\sigma-\tau-|x-y|)}^
  {L(\rho-\sigma-\tau+|x-y|)} \dfrac{d}{ds}\widehat{\chi'}(s)\,ds 
  \, \Big|\,d\rho\ \le\  2\bignorm[\big(\widehat{\chi'}\big)'][1]
\end{equation*}
by taking the absolute value inside again and using Fubini's theorem.
This bound is independent of all other variables, including $L$, so by 
Theorem~\ref{thm:Fourier} the size of \eqref{eq:boundary} is controlled
by $|t|^{-3/2}\norm[f][1]\norm[g][1]$.  The proof of 
Theorem~\ref{thm:dispersive} is complete.

{\bf Remark.}
The second part of \eqref{eq:IBP} gives the impression of being a boundary
term, so it would be satisfying to see it vanish as $L \to \infty$.
An additional estimate for \eqref{eq:boundary} shows that this occurs.
After absolute values are brought inside, the integral in $\rho$ can also
be bounded above by $2\norm[\widehat{\chi'}][1]/(L|x-y|)$ by assuming no
cancellation between the evaluations of $\widehat{\chi'}$.  This
provides pointwise (in $(x,\sigma,y,\tau)$) convegence to zero as
$L\to\infty$, and the bound used above lets us apply dominated convergence.

\bibliographystyle{amsplain}

\medskip\noindent
\textsc{Division of Astronomy, Mathematics, and Physics, 253-37 Caltech, Pasadena, CA 91125, U.S.A.}\\
{\em email: }\textsf{\bf mikeg@its.caltech.edu} 

\end{document}